\documentclass[12pt]{article}
\usepackage{amsmath}
\usepackage{amsfonts}
\newtheorem{theorem}{Theorem}[section]
\newtheorem{lemma}{Lemma}[section] 
\newtheorem{corollary}{Corollary}[section] 
\newtheorem{proposition}{Proposition}[section]

\newcommand{\ts}[1]{\langle #1\rangle}

\title{\bf On locally solvable   maximal subgroups of the multiplicative group of a division ring }

\author{Bui Xuan Hai 
\\{\small\em Faculty  of Mathematics and Computer Science, University of  Science}\\{\small\em VNU - HCM-City}\\ {\small\em 227 Nguyen Van Cu str., Dist. 5, Ho Chi Minh City, Vietnam}\\ {\small\em e-mail:bxhai@hcmus.edu.vn}\\Dang Vu Phuong Ha\\{\small \em Department of Basic Sciences, University of Architecture}\\{\small \em 196 Pasteur Str., Dist. 1, HCM-City, Vietnam}\\{\small \em e-mail: phuongha057@yahoo.com}}

\date{}
\rm\em 
\normalsize
\begin{document}
\maketitle
\newcommand{\dpcm}{ \hfill \rule{3mm}{3mm}}
\def\Box{\dpcm}
\newsymbol \lneq  2308   
\newsymbol \lsubset 2320    
\def\xd{\linebreak} 
\begin{abstract} 
Let $D$ be  a division ring  and $D^*$ be the multiplicative group of $D$. In this paper we study locally solvable maximal subgroups of $D^*$. 
\end{abstract}



\newpage

\section{Introduction and convention}

Let $D$ be a division ring with the center $F$. All subgroups  considered in this paper are subgroups of the multiplicative group $D^*$ of $D$. So,  sometime we say that $G$ is a subgroup of $D$ with the understanding that $G$ is in fact a subgroup of $D^*$. One of the well-known  results of L.K. Hua states that if $D^*$ is solvable, then $D$ is commutative. Here, first we generalize this classical result by proving that if $D^*$ is locally solvable, then $D$ is also commutative. Using this fact, we obtain a series of other results, concerning locally solvable subgroups of $D^*$. The paragraph 2 is devoted to the study of properties of locally solvable maximal subgroups in $D$. Among the new obtained results  we would like to notice  the following fact: {\em  Let $D$ be a non-commutative division ring with the center $F$ and suppose that $M$ is  a  non-abelian locally solvable maximal subgroup of $D^*$. If $M'$ is algebraic over $F$, then $[D: F] <\infty$.}  Before, this fact was proved in [3, Th. 6]
with a stronger supposition of the algebraicity of $M$. Here, we replace the condition of algebraicity of $M$ by the algebraicity of derived subgroup $M'$ of $M$. 

Finally, recall that in [1] it was proved that in a non-commutative centrally finite division ring $D$, every nilpotent maximal subgroup of $D^*$ is the multiplicative group of some maximal subfield of $D$. More generally, in [4, Cor. 5], the author proved that in an arbitrary division ring, every nilpotent maximal subgroup is abelian. On the other hand, it is easy to see that every abelian maximal subgroup is the multiplicative group of  some  maximal subfield of  $D$ (see Pro. 2.2 in the text).  In [5, Th. 3.2] this result was carried over for  locally nilpotent maximal subgroups of  a non-commutative division ring $D$ that is algebraic over its center $F$. Now, using the result, mentioned above, here   we  can  replace the algebraicity of $D$ by the algebraicity of $M'$, but the obtained result  is the same. In fact, we shall show that if $M$ is a locally nilpotent maximal subgroup such that $M'$ is algebraic over $F$, then $M$ is  the multiplicative group of some maximal subfield of $D$.  

Throughout this paper the following notations will be used consistently: For a division ring $D$ we  denote by $D^*$  and $F$ its multiplicative group  and  center respectively.  We say that a division ring $D$ is {\em algebraic over  $F$} if every element of $D$ is algebraic over $F$.  A division ring $D$ is called {\em centrally finite} provided it is a finite dimensional vector space over $F$. If $R$ is a ring with identity $1\neq 0$  and $A$ is a  nonempty subset of $R$, then $ C_R(A)$ denotes the centralizer of $A$ in $R$, i.e.
$$C_R(A):=\{x\in R\vert~ xa=ax \mbox{ for all } a\in A\}.$$
If $G$ is a group, then $Z(G)$ is the center  and $G'$ is the drived subgroup of $G$ respectively.   For a subgroup $G$ of $D^*$, denote by $F[G]$ and $F(G)$ the subring and the division subring respectively of $D$ generated by the set $F\cup G$. We say that a subgroup $G$ of $D^*$  is {\em irreducible} (resp. {\em absolutely irreducible}) if $F(G)=D$ (resp. $F[G]=D$).  All  another notation and symbols in this paper are standard and one can find, for example, in [6], [7] and [8]. 

\section{Locally solvable  maximal subgroups}

Our purpose is to study  the properties of locally solvable maximal subgroups in a division ring $D$. For the convenience we restate  the following two results of Wehrfritz which will be used in this work.\\

\noindent
{\bf Theorem A.} [11, Cor. 4] {\em Let  $H$ be a locally solvable normal subgroup of the  absolutely irreducible  subgroup $G$ of $GL_n(D)$.  If either $n=1$ or $H=G$, then $H$ is abelian-by-locally finite and $G/C_G(H)$ is abelian-by-periodic.} \\

\noindent
{\bf Theorem B.} ([10, Th. 5.7.11, p. 215]) {\em Let $H$ be a locally nilpotent normal subgroup of the absolutely irreducible subgroup $G$ of $GL_n(D)$. Then $H$ is centre by locally-finite and $G/C_G(H)$ is periodic.}\\

We note also the following two  very simple lemmas we need. 

\begin{lemma}\label{lem 2.1} Let $D$ be a division ring with the center $F$ and $G$ be a subgroup of $D^*$. If $N$ is a normal subgroup of $G$ and $L=F(N)$, then $G$ is contained in $N_{D^*}(L^*)$. In particular, if $G=D^*$, then $L$ is normal in $D^*$.
\end{lemma}

 \noindent
{\em  Proof.} Suppose that $x\in G$. Since $N\triangleleft G, xNx^{-1}=N\subseteq L$, so $N\subseteq x^{-1}Nx$. Since it is obvious that $F\subseteq xLx^{-1}$, we have $xLx^{-1}=L$. So, it follows that $G\leq N_{D^*}(L^*)$.\dpcm

\begin{lemma}\label{lem 2.2}  Let $D$ be a division ring with the center $F$. If   $G$ is  an irreducible  subgroup of $D^*$, then $C_D(G)=F$. 
\end{lemma}

\noindent
{\em  Proof.} It is obviuos that  $C_D(G)=C_D(F(G))=C_D(D)=F$.\dpcm\\

The following theorem generalizes the classical result of Hua, mentioned in the Introduction  (see, for example [6, p. 223] or also the generalization of this result by Stuth in  [9]).

\begin{theorem}\label{th 2.1} If $D^*$ is locally solvable, then $D$ is a field.
\end{theorem}

\noindent
{\em  Proof.} By applying of Theorem A for $n=1, H=G=D^*$, we can find  an abelian normal subgroup $A$ of $D^*$ such that $D^*/A$ is locally finite. By Lemma 2.1, $F(A)^* \triangleleft D^*$. In view of  Cartan-Brauer-Hua  Theorem (see, for example [6, (13.17), p. 222]), either $F(A)=D$ or $F(A)\subseteq F$. If $F(A)=D$, then by Lemma 2.2, we have $C_D(A)=F$. Since $A$ is abelian, it follows that $A\subseteq F$. So, in both cases we have $A\subseteq F$. Since $D^*/A$ is locally finite and $A\subseteq F$, it follows that $D^*/F^*$ is locally finite. In particular, $D$ is radical over $F$. Now, in virtue of  Kaplansky's Theorem (see, for example [6, (15.15), p. 259]), $D$ is a field.\dpcm

\begin{proposition} \label{pro 2.1} Let $D$ be a non-commutative division ring and suppose that $M$ is a locally solvable  maximal subgroup of $D^*$. Then, the  following statements hold:

(i) $F^*\lneq M$.

(ii) If $A$ is an  abelian normal subgroup of $M$ and $L=F(A)$, then $L$ is a subfield of $D$ and $L^*\triangleleft M$.
\end{proposition}

\noindent
{\em  Proof.}  (i) Since $M$ is maximal in $D^*$, either $F^*M=M$ or $F^*M=D$.  If $F^*M=D$, then $D^*$ is locally solvable. So, by Theorem 2.1, $D$ is commutative, that contradicts to the supposition. Hence,  $F^*M=M$ and consequently   $F^*\leq M$. Since $D$ is non-commutative and $M$ is maximal in $D^*$, $F^*\neq M$. 

(ii) By Lemma 2.1, $M\leq N_{D^*}(L^*)$. So, to prove $L^*\triangleleft M$, it suffices to show that $L^*\leq M$. Since $M$ is maximal in $D^*$ and $M\leq N_{D^*}(L^*)\leq D^*$, either $N_{D^*}(L^*)=D^*$ or $N_{D^*}(L^*)=M$. If $N_{D^*}(L^*)=D^*$, then by Cartan-Brauer-Hua Theorem , either $L\subseteq F$ or $L=D$. If $L=D$, then by Lemma 2.2, $C_D(A)=F$ and it follows that $A\subseteq F$ and consequently $L^*=F(A)^*=F^*\leq M$. Now, if $N_{D^*}(L^*)=M$, then we also have $L^*\leq M$. Thus, $L^*\leq M$ in any case, as we desired to show. Since $M$ is locally solvable, $L^*$ is locally solvable too. Hence, in view of  Theorem 2.1, $L$ is a field, as we desired to prove.\dpcm

\begin{proposition}\label{pro 2.2}  Let $D$ is a non-commutative division ring with the center $F$ and suppose that $M$ is a maximal subgroup of $D^*$. If $M$ is abelian, then $M$ is the  multiplicative  group of some maximal subfield of $D$.
\end{proposition}

\noindent
{\em Proof.} By Proposition 2.1, $F^*\leq M$. By [1, Pro. 1], either $F(M)=D$ or $M\cup\{0\}$ is a division subring of $D$. Since $D$ is non-commutative and $M$ is abelian, the first assertion can not occur. So, $M\cup\{0\}$ is a division subring of $D$. Therefore, $M\cup\{0\}$ is a maximal subfield of $D$. \dpcm

\begin{lemma}\label{lem 2.3} Let $D$ be a non-commutative division ring with the center $F$ and suppose that $M$ is a locally solvable maximal subgroup of $D^*$. If $M$ is irreducible, then $Z(M)=F^*$.
\end{lemma}

\noindent
{\em Proof.} By Proposition 2.1, $F^*\leq M$, so it follows that $F^*\leq Z(M)$. Since  $M$ is irreducible, $F(M)=D$. Therefore  $Z(M)\leq F^*$ and consequently, $Z(M)=F^*$.\dpcm\\

We note  the following simple fact, whose proof will be omitted.

\begin{lemma}\label{lem 2.4} Let $M$ be a non-abelian group and suppose that $A$ is an abelian normal subgroup of $M$. Then $A$ is contained in some maximal abelian normal  subgroup of $M$. In particular, $M$ contains maximal abelian normal  subgroups.
\end{lemma}

\begin{proposition}\label{pro 2.3} Let $D$ be a non-commutative division ring with the center $F$ and suppose that $M$ is a locally solvable maximal subgroup of $D^*$. Then, $M$ is irreducible if and only if $M$ is non-abelian.
\end{proposition}

\noindent
{\em Proof.} Suppose that $M$ is irreducible. By Lemma 2.3, $Z(M)=F^*$. If $M$ is abelian, then $M=F^*$ and it follows that $D=F(M)=F$, that is a contradiction. 

Conversely, suppose that $M$ is non-abelian. If $M$ is not irreducible, then $F(M)^*=M$. Since $M$ is locally solvable, by Theorem 2.1, $F(M)$ is a field and consequently $M$ is abelian that is a contradiction.\dpcm 

\begin{theorem}\label{th 2.2} Let $D$ be a non-commutative  division ring with the center $F$ and suppose that $M$ is a  non-abelian locally solvable maximal subgroup of $D^*$ and $A$ is a maximal  abelian normal  subgroup of $M$. By setting $L=F(A), K=C_D(L), V=C_M(A)$, we have:

(i) $A=L^*$.

(ii) Either $A=F^*$ or $V=L^*=K^*$ and $K$ is a maximal subfield of $D$.
\end{theorem}

\noindent
{\em Proof.} First, note that in view of Lemma 2.4, such a subgroup $A$ of $M$ exists.  

(i) By Proposition 2.1, $A\triangleleft L^*\triangleleft M$. Since $A$ is  maximal  abelian normal in $M$, it follows that $L^*=A$. 

(ii) By Lemma 2.1, $M\leq N_{D^*}(L^*)\leq D^*$. In view of the maximality of  $M$  in $D^*$, either $N_{D^*}(L^*)=D^*$ or $N_{D^*}(L^*)=M$. If $N_{D^*}(L^*)=D^*$, then as in the proof of Proposition 2.1 we have $L=F$. So, by (i), $A=F^*$. Now, suppose that $N_{D^*}(L^*)=M$. Since $K=C_D(L)$, we have $K^*\leq N_{D^*}(L^*)=M$. Therefore $K^*$ is locally solvable, so by Theorem 2.1, $K$ is a field. Thus, we  have $A=L^*\triangleleft K^*\leq M$. Suppose that $K_1$ is a subfield of $D$ containing $K$. Since $L\subseteq K\subseteq K_1, K_1\subseteq C_D(L)=K$. Hence $K_1=K$ and we conclude that $K$ is a maximal subfield of $D$. Since $K^*\leq M$ and $A=L^*\leq K^*$, it follows that $K^*=C_M(A)=V$. Suppose that $x\in M, y\in K^*$ and $a\in A$ are arbitrary. Then, we have
$$x^{-1}yxa=x^{-1}yxax^{-1}x=x^{-1}xax^{-1}yx=ax^{-1}yx.$$

This shows that $x^{-1}yx\in C_M(A)=K^*$. So, $K^*\triangleleft M$. By  maximality of $A$, this forces $K^*=A$. Thus, $V=L^*=K^*$ as it was required to prove.\dpcm

\begin{proposition}\label{pro 2.4} Let $D$ be a non-commutative  division ring with the center $F$. If   $M$ is a non-abelian  locally solvable maximal subgroup of $D^*$, then $M$ contains a unique maximal abelian normal subgroup. 
\end{proposition}

\noindent
{\em Proof.} By  Lemma 2.4, $M$ contains maximal abelian normal  subgroups. Now, suppose that $A$ is an arbitrary maximal abelian normal  subgroup of $M$. By Theorem 2.2 (i), $F^*\leq A$. If $A=F^*$ and $B$ is an another maximal abelian normal subgroup of $M$, then $F^*=A\leq B\triangleleft M$; hence $B=A=F^*$. Therefore, in this case $F^*$ is a unique maximal abelian normal subgroup of $M$. Thus, we can suppose that $A\neq F^*$ for every maximal abelian normal subgroup $A$ of $M$. Now, suppose that $A_1, A_2$ are maximal abelian normal  subgroups of $M$. By Theorem 2.2,  $K_1=A_1\cup\{0\}$ and $K_2=A_2\cup\{0\}$ are maximal subfields of $D$. 

Suppose that $x\in K_1^*\setminus K_2^*$. Then, there exists some element $y\in A_2$ such that $xy\neq yx$ or equivalently, $a=xyx^{-1}y^{-1}\neq 1$. For any $b\in A_2$, since elements $y^{-1}, xyx^{-1}$ are both in $A_2$, we have
$$ab=xyx^{-1}y^{-1}b=xyx^{-1}by^{-1}=bxyx^{-1}y^{-1}=ba.$$

Hence $a\in K_2$. Since $x+1\in K_1^*$ and $(x+1)y\neq y(x+1)$, by the similar way as above we can conclude that $1\neq c=(x+1)y(x+1)^{-1}y^{-1}\in K_2^*$. Then, we have 

\hspace*{5cm}$cy(x+1)=(x+1)y$

\hspace*{4.1cm}$\Longleftrightarrow cyx+cy=xy+y$

\hspace*{4.1cm}$\Longleftrightarrow cyx-ayx=(1-c)y$

\hspace*{4.1cm}$\Longleftrightarrow y(c-a)x=y(1-c)$

\hspace*{4.1cm}$\Longleftrightarrow (c-a)x=1-c.$

If $c\neq a$, then $x\in A_2=K_2^*$. If $c=a$, then $a=c=1$. Thus, in both cases we  have  a contradiction. Therefore $K_1^*\subseteq K_2^*$. By symmetry, $K_2^*\subseteq K_1^*$, so $K_2^*= K_1^*$ or $A_1=A_2$.\dpcm

\begin{proposition}\label{pro 2.5} Let $D$ be a non-commutative division ring with the center $F$ and suppose that $M$ is a non-abelian  locally solvable maximal subgroup of $D^*$. Then, $D$ is centrally finite if and only if there exists some maximal subfield $K$ of $D$ such that $F^*\lneq K^*\triangleleft M$ and $[M: K^*] < \infty$.
\end{proposition}

\noindent
{\em Proof.} The ``if"  follows from [3, Th. 6].  For the ``only if", note that, since $F^*\lneq K^*\triangleleft M$ and $K$ is a field, it is easy to see that $F(K^*)$ is a subfield which is contained in $M$. Since $[M: K^*] <\infty, [M: F(K^*)^*] < \infty$. By [3, Lemma 6], $[D: F] < \infty$.\dpcm

\begin{proposition}\label{pro 2.6} Let $D$ be a non-commutative centrally finite division ring with the center $F$ and  suppose that $M$ is a   locally solvable maximal subgroup of $D^*$. If $M$ is non-abelian, then $M$ is not radical over $F$.
\end{proposition}

\noindent
{\em Proof.} Suppose that $M$ is non-abelian locally solvable maximal subgroup of $D^*$. By [3, Th. 6], there exists some maximal subfield $K$ of $D$ such that $F^*\lneq K^*\leq M$ and $K/F$ is a Galois extension. If $M$ is radical over $F$, then $K$ is radical over $F$ too. By [6, (15.13), p. 258], the prime subfield $P$ of $F$ has the  characteristic $p > 0$ and either $K$ is purely inseparable over $F$ or $K$ is algebraic over $P$. Since $K/F$ is Galois, the first case cannot occur, so $K$ is algebraic over $P$. Consequently, $D$ is algebraic over the  finite field $P$ and by well-known theorem of Jacobson [see 6, (13.11), p. 219], $D$ is a field, that is a contradiction.\dpcm

\begin{lemma}\label{lem 2.5} Let $D$ be a non-commutative division ring with the center $F$ and suppose that $M$ is a  locally solvable maximal subgroup of $D^*$ and $x\in M\setminus F$. If $M'\subseteq F$, then $F(x)^*\triangleleft M$.
\end{lemma}

\noindent
{\em Proof.} By Proposition 2.1, $F^*\lneq M$, so such an  $x$ exists. For any $m\in M$ we have $m^{-1}xm=xx^{-1}m^{-1}xm=x[x, m]\in xM'\subseteq F(x)$, so $x\in mF(x)m^{-1}$. Hence $m^{-1}F(x)m\subseteq F(x), \forall m\in M$ and consequently  $M\leq N_{D^*}(F(x)^*)\leq D^*$. Since $M$ is maximal in $D^*$, either $ N_{D^*}(F(x)^*)=M$ or $N_{D^*}(F(x)^*)=D^*$. If $N_{D^*}(F(x)^*)=D^*$, then by Cartan-Brauer-Hua Theorem, either $F(x)=D$ or $F(x)\subseteq F$. Since  these cases are  both impossible, $N_{D^*}(F(x)^*)=M$ or $F(x)^*\triangleleft M$.\dpcm

\begin{theorem}\label{th 2.3} Let $D$ be a non-commutative division ring with the center $F$ and suppose that $M$ is an   irreducible   locally solvable maximal subgroup of $D^*$.  If $M$ is metabelian, then the following statements hold:

(i) $M$ contains a unique  maximal  abelian normal subgroup $A$ such that $F^*\lneq A$ and $M'\leq A$.

(ii) $K=A\cup \{0\}$ is a maximal subfield of $D$. 
\end{theorem}

\noindent
{\em Proof.} (i) In view of Proposition 2.3, $M$  is   non-abelian.  By Proposition 2.4 and Theorem 2.2,  $M$ has a unique  maximal abelian normal subgroup containing $F^*$, say  $A$.  Since $M$ is metabelian, $A$ must contain $M'$. Suppose that 
 $A=F^*$. Then $M'\leq F^*\lneq M$. For $x\in M\setminus F$, by Lemma 2.5, we have $F(x)^*\triangleleft M$, so by Theorem 2.1, $F(x)$ is a field. Since  $A$ is maximal abelian normal in $M$, we have  $F(x)^*=F^*$ or $x\in F$, that is a contradiction. 
Therefore $A\neq F^*$.  

(ii) By Theorem 2.2, $K$ is a maximal subfield of $D$. \dpcm

\begin{lemma}\label{lem 2.6} Let $D$ be a non-commutative division ring with the center $F$. If $M$ is an irreducible locally solvable maximal subgroup of $D^*$, then $M'\not\subseteq F$.
\end{lemma}

\noindent
{\em Proof.}   As it was shown in the proof of Theorem 2.3, $M$ is non-abelian and  we have a unique abelian normal subgroup $A$ of $M$ containing $F^*$. Suppose that $M'\subseteq F$ and $x\in M\setminus F$. By Lemma 2.5, $F(x)^*\triangleleft M$. In view of  Theorem 2.1, $F(x)$ is a field and consequently $F(x)^*\leq A $. Therefore, $x\in A$ for any $x\in M\setminus F$. Since $F^*\leq A$ , it follows that $M= A$. In particular, $M$ is abelian, that is a contradiction.\dpcm

\begin{proposition}\label{pro 2.7} Let $D$ be a non-commutative division ring with the center $F$. Suppose that $M$ is an  irreducible  metabelian locally solvable maximal subgroup of $D^*$ and $A$ is a maximal abelian normal  subgroup of $M$. If a subgroup $N$ of $M$ strictly contains $A$, then $N$ is irreducible.
\end{proposition}

\noindent
{\em Proof.} By Theorem 2.3, $K=A\cup\{0\}$ is a maximal subfield of $D$ and $M'\subseteq K, F\subseteq K$. Therefore  $M'\leq K^*\lneq N$. Clearly $N\triangleleft M$, so by Lemma 2.1, 
$$M\leq N_{D^*}(F(N)^*)\leq D^*.$$

Suppose that $M= N_{D^*}(F(N)^*)$. Then, in view of  Theorem 2.1, $F(N)$ is a subfield of $D$ containing $K$. Since $K$ is a maximal subfield of $D$, it follows that $K=F(N)$. However, this is  impossible since $N\neq K$.  Hence $D=N_{D^*}(F(N)^*)$. By Cartan-Brauer-Hua Theorem, either $F(N)\subseteq F$  or $F(N)=D$. Since $ F(N)\neq F, F(N)=D$ or $N$ is irreducible.\dpcm\\

If $K$ is a subfield of a division ring $D$, then we denote by $[D: K]_l$ ($[D: K]_r$, resp.) the dimension of left (right, resp.) vector space $D$ over $K$. 

\begin{lemma}\label{lem 2.7} Let $D$ be a division ring with the center $F$ and suppose that $K$ is a subfield of $D$ containing $F$. If $[D: K]_l < \infty$ or $[D: K]_r < \infty$, then $[D: F] < \infty$.
\end{lemma}

\noindent
{\em Proof.} Consider the case, when  $[D: K]_l < \infty$. By putting  $L=C_D(K)$, we have $K\subseteq L$. Since $[D: K]_l < \infty, [D: L]_l < \infty$. By [6, (15.4), p. 253], $[K: F] <\infty$. Hence $[D: F] < \infty$. The remaining  case may be  considered similarly. \dpcm

\begin{theorem}\label{th 2.4} Let $D$ be a non-commutative division ring with the center $F$. Suppose that $M$ is an   irreducible  metabelian locally solvable maximal subgroup of $D^*$ and $K$ is a   maximal subfield of $D$.  If there exists some algebraic over $K$ element $\alpha\in M\setminus K$, then $D$ is a finite dimensional vector space over $F$.
\end{theorem}

\noindent
{\em Proof.} By Theorem 2.3, $F^*\lneq K^*\triangleleft M$. Consider the minimal polynomial of $\alpha$ over $K$
$$f(X)=X^n+b_{n-1}X^{n-1}+\ldots+b_1X+b_0\in K[X]$$
and put $R:=\sum_{i=0}^{n-1}K\alpha^i$. Clearly $R$ is a left vector space over $K$ of dimension $n$ with the  basis $\{1, \alpha, \alpha^2, \ldots, \alpha^{n-1}\}$. Since $K^*\triangleleft M$ and $\alpha\in M$, it is easy to see that $R$ is a subring of $D$. Now, suppose that $0\neq x\in R$. Then, there exists some positive integer $m$ such that
$$1=c_1x+c_2x^2+\ldots+c_mx^m \mbox{ for } c_1, c_2, \ldots, c_m\in K.$$

Therefore $1=(c_1+c_2x+\ldots+c_mx^{m-1})x$, so $x$ is invertible. Thus $R$ is a division subring of $D$. Let $N=\ts{K^*, \alpha}$ be the subgroup of $M$ generated by $K^*$ and $\alpha$. Since $\alpha\not\in K, N$ strictly contains $K^*$. By Proposition 2.7, $F(N)=D$.  On the other hand, $R=F(N)$, hence $D=R$. So $[D: K]_l <\infty$. Now, by Lemma 2.7,  $[D: F] < \infty$.\dpcm

\begin{theorem}\label{th 2.5} Let $D$ be a non-commutative division ring with the center $F$. Suppose that $M$ is a non-abelian  locally solvable maximal subgroup of $D^*$ and $A$ is a maximal abelian normal subgroup of $M$. If there exists some element $\alpha\in A\setminus F$ such that $\alpha$ is algebraic over $F$, then $[D: F] < \infty$.
\end{theorem}

\noindent
{\em Proof.} By Theorem 2.2, we have $F^*\leq A=F(A)^*\triangleleft M$. Denote by $f(X)=min(F, \alpha)$ the minimal polynomial of $\alpha$ over $F$ and suppose that
$$f(X)=X^n+a_{n-1}X^{n-1}+\ldots+a_1X+a_0.$$

For any $b\in M$ we have $b^{-1}\alpha b\in A$ and $f(b^{-1}\alpha b)=b^{-1}f(\alpha)b=0.$ So, by applying [8, 3.3.3, p. 53], it follows that
$$[M: C_M(\alpha)]=|\alpha^M |< deg(f)=n < \infty.$$ 

Set $N_1=Core(C_M(\alpha))$. By [8, 3.3.5, p. 53] we have
$$N_1\triangleleft  M, N_1\leq C_M(\alpha) \mbox{ and } [M: N_1] < \infty.$$

By Lemma 2.2, $M\leq N_{D^*}(F(N_1)^*)\leq D^*$. Hence, either $N_{D^*}(F(N_1)^*)=D^*$ or $N_{D^*}(F(N_1)^*)=M$.
If  $N_{D^*}(F(N_1)^*)=D^*$, then by Cartan-Brauer-Hua Theorem, either $F(N_1)\subseteq F$ or $F(N_1)=D$. Suppose that $F(N_1)=D$. Then, by Lemma 2.2, $F=C_D(N_1)$. Since $N_1\leq C_M(\alpha), \alpha \in C_D(N_1)$. So, $\alpha \in F$, that is a contradiction. Thus, $F(N_1)\subseteq F$, so $N_1\subseteq F$. Since $[M: N_1] < \infty, [M: F^*] < \infty$. By [1, Cor. 4, p. 426], $D$ is commutative, that contradicts to the supposition. Thus, $N_{D^*}(F(N_1)^*)=M$. Then, $F(N_1)^*\leq M$, so $F(N_1)^*$ is locally solvable and by Theorem 2.1, $F(N_1)$ is a field. Now, since $[M: F(N_1)^*] < \infty$, in view of [3, Lemma 6], $[D: F] < \infty$.\dpcm

\begin{theorem}\label{th 2.6} Let $D$ be a non-commutative division ring with the center $F$. Suppose that $M$ is an   irreducible  metabelian  locally solvable maximal subgroup of $D^*$. If there exists some algebraic over $F$ element $\alpha\in M\setminus F$, then $D$ is a finite dimensional vector space over $F$.
\end{theorem}

\noindent
{\em Proof.} By Theorem 2.3, there exists some maximal subfield $K$ of $D$ such that $F^*\leq K^*\triangleleft M$ and $A=K^*$ is a maximal abelian normal  subgroup of $M$. Suppose that $\alpha\in M\setminus F$. If $\alpha\not\in K^*$, then by Theorem 2.4 , $[D: F] < \infty$. If $\alpha\in K$, then by Theorem 2.5, $[D: F] <\infty$.\dpcm

\begin{lemma}\label{lem 2.8} Let $D$ be a non-commutative division ring with the center $F$ and suppose that $M$ is an irreducible locally solvable maximal subgroup of $D^*$. If $M'$ is algebraic over $F$, then $M$ is absolutely irreducible.
\end{lemma}

\noindent
{\em Proof.} Recall that $M$ is non-abelian in view of  its  irreducibility.  Consider   the  group $F[M]^*$of all units of the ring $F[M]$. Since $M$ is maximal in $D^*$, it follows that either $F[M]^*=M$ or $F[M]^*=D^*$. Suppose that $F[M]^*=M$. For any $a, b\in M'$, since $a, b^{-1}$ are algebraic over $F, a\pm b, ab^{-1}$ are algebraic over $F$ too. Therefore $F(ab^{-1})=F[ab^{-1}]\subseteq F[M]^*=M$. It follows that $ab^{-1}\pm 1\in M$ and consequently $a\pm b\in M$. Hence, $a\pm b\pm c\in M$ for any $a, b, c\in M'$. By Proposition 2.1,  $F^*\leq M$, so we can conclude that $F(M')^*\leq M$. Since $M$ is locally solvable, $F(M')^*$ is locally solvable too. Therefore, by Theorem 2.1, $F(M')$ is a field. In particular, $M'$ is abelian and so $M$ is metabelian. By Lemma 2.6 , $M'\not\subseteq F$. Hence, there exists some element $x\in M'\setminus F$ such that $x$ is algebraic over $F$. By Theorem  2.6  , $[D: F] <\infty$, so $F[M]=F(M)$ by [5, Lem. 2.3]. By supposition $F[M]^*=M$, so we have $F(M)^*=M$. By Theorem 2.1, $M$ is abelian, that is a contradiction. Thus, $F[M]^*=D^*$ or $F[M]=D$, i.e. $M$ is absolutely irreducible as it was required to prove.\dpcm

\begin{corollary}\label{cor 2.1} Let $D$ be a non-commutative division ring with the center $F$ and suppose that $M$ is an irreducible locally solvable maximal subgroup of $D^*$. If $M'$ is algebraic over $F$, then $D$ is a finite dimensional vector space over $F$.
\end{corollary}

\noindent
{\em Proof.} By Lemma 2.8, $M$ is absolutely irreducible. By Theorem A, there exists an   abelian normal subgroup $N$ of $M$ such that $M/N$ is locally finite. Denote by $A$ a maximal  abelian normal subgroup of $M$ containing $N$. Then, $M/A$ is locally finite. If $A=F^*$, then by [3, Th. 6], $[D: F] < \infty$. If $A\neq F^*$, then in view of  Theorem 2.3, $M'\leq A$ and by Lemma 2.6, $M'\not\subseteq F$. Hence, by Theorem 2.6, $[D: F] < \infty.$\dpcm

\begin{theorem}\label{th 2.7} Let $D$ be a non-commutative division ring with the center $F$ and suppose that $M$ is an   irreducible  maximal subgroup of $D^*$ such that $M'$ is algebraic over $F$. If  $M$ is locally solvable, then there exists some maximal subfield $K$ of $D$ such that the following conditions hold:

(i) $K^*\triangleleft M$.

(ii) $K$ is a Galois extension of $F$.

(iii) $M/K^*\simeq Gal(K/F)\simeq \mathbb{Z}_p$, where $p$ is a prime number.
\end{theorem}

\noindent
{\em Proof.} Suppose that $M$ is an irreducible  maximal subgroup of $D^*$ and $M'$ is algebraic over $F$. If $M$ is locally solvable, then by Corolarry 2.1, $[D: F] < \infty$. Hence, by [3, Th. 6], there exists some maximal subfield $K$ of $D$ satisfying the conditions (i), (ii) and (iii).\dpcm

\begin{theorem}\label{th 2.8} Let $D$ be a non-commutative division ring with the center $F$ and suppose that $M$ is an  irreducible  locally solvable maximal subgroup of $D^*$. Then, $M'$ is algebraic over $F$ if and only if $M$ is metabelian and there exists some element $x\in M\setminus F$ such that $x$ is algebraic over $F$.
\end{theorem}

\noindent
{\em Proof.} Suppose that $M'$ is algebraic over $F$. By Theorem 2.7, there exists some maximal subfield $K$ of $D$ such that $K^*\triangleleft M$ and $M/K^*\simeq \mathbb{Z}_p$. Then, $M'\leq K^*$, so $M$ is metabelian. By Lemma 2.6, there exists some element $x\in M'\setminus F$ and $x$ is algebraic over $F$ by supposition.

Conversely, suppose that  $M$ is metabelian and  there exists some element $x\in M\setminus F$ such that $x$ is algebraic over $F$. By Theorem 2.6, $[D: F] < \infty$. In particular, $M'$ is algebraic over $F$.\dpcm

\begin{lemma}\label{lem 2.9} Let $D$ be a non-commutative division ring with the center $F$ and suppose that $M$ is a maximal subgroup of $D^*$ containing some  non-central element $a$ algebraic over $F$. If $M$ is abelian, then $[D: F] < \infty.$
\end{lemma}

\noindent
{\em Proof.} Since $M$ is abelian and $D$ is non-commutative, by applying Proposition 2.2 we  conclude that   $F(M)=M\cup\{0\}$ is the maximal subfield of $D$. Putting  $L=C_D(F(a))$, we have $M\leq L^*\leq D^*$. Then, either $L^*=D^*$ or $L^*=M$. If $L^*=D^*$, then $D=C_D(F(a))$, so $a\in F$, that  contradicts to the supposition. Hence $L^*=M$, so by [6, (15.7), p. 254], $C_D(L)=L$. Since $C_D(L)=C_D(C_D(F(a))$ and $[F(a): F] < \infty$, by Double Centralizer Theorem we have  $L=C_D(L)=F(a)$. By applying [6, (15.4), p. 253], $[D: F] < \infty.$\dpcm\\

In [3, Th. 6], it was proved that if $M$ is a non-abelian locally solvable maximal subgroup of $D^*$ and $M$ is algebraic over the center $F$ of $D$, then $D$ is a finite dimensional vector space over $F$. Now, we are ready to  show  that this result remains   also true if we replace the condition of algebraicity of $M$ by the algebraicity of derived subgroup $M'$ of $M$. 

\begin{theorem}\label{th 2.9} Let $D$ be a non-commutative division ring with the center $F$ and suppose that $M$ is  a  non-abelian locally solvable maximal subgroup of $D^*$. If $M'$ is algebraic over $F$, then $[D: F] <\infty$. 
\end{theorem}

\noindent
{\em Proof.} If $M$ is irreducible, then by Corollary 2.1, $[D: F] <\infty$. If $M$ is not irreducible, then $F(M)^*=M$ and by Theorem 2.1, $F(M)$ is a field. Hence, by Lemma 2.9, $[D: F] < \infty.$\dpcm 

\begin{corollary}\label{cor 2.2}  Let $D$ be a non-commutative division ring with the center $F$. If $D^*$ contains some locally solvable maximal subgroup $M$ such that $M'$ is algebraic over $F$, then every locally solvable subgroup of $D^*$ is solvable.
\end{corollary}

\noindent
{\em Proof.} By Theorem 2.9, $D$ is centrally finite, so $D^*$ and its subgroups can be considered as a linear groups over $F$. Since  every locally solvable linear group is solvable, the proof is now complete. \dpcm\\

 In [5, Th. 3.2] it was proved that, every locally nilpotent maximal subgroup of $D^*$ is the multiplicative group of some maximal subfield of $D$ provided  $D$  is algebraic over its center $F$. Now, using the results obtained above,  we shall  show that this result remains also true with the  weaker condition for the algebraicity of derived subgroup $M'$ of $M$ instead of the algebraicity of $D$. 

\begin{theorem}\label{th 3.1} Let $D$ be a division ring with the center $F$ and suppose that $M$ is a locally nilpotent maximal subgroup of $D^*$. If $M'$ is algebraic over $F$, then $M$ is the multiplicative group of some maximal subfield of $D$.
\end{theorem}

\noindent
{\em Proof.} By Proposition 2.2, it suffices to prove that $M$ is abelian. Thus, suppose that $M$ is non-abelian. Then, by Proposition 2.3, $M$ is irreducible. Hence by Lemma 2.8, $M$ is absolutely irreducible. By applying Theorem B for $H=G=M$, we conclude that $M/Z(M)$ is torsion group. On the other hand, by Lemma 2.3, $Z(M)=F^*$, so $M/F^*$ is torsion. Therefore $M$ is radical over $F$. Since $M'$ is algebraic over $F$, by Theorem 2.9, $[D: F] < \infty$. By Proposition 2.6, $M$ is not radical over $F$, that contradicts to the conclusion above. Hence $M$ is abelian, as we desired to prove.\dpcm\\

\end{document}